\documentclass[11pt]{article}

\usepackage[T1]{fontenc}
\usepackage{lmodern}
\usepackage[margin=1in]{geometry}
\usepackage{amsmath,amssymb}
\usepackage{booktabs}
\usepackage{float}
\usepackage{algorithm}
\usepackage{algpseudocode}
\usepackage[numbers,sort&compress]{natbib}
\usepackage{microtype}
\usepackage{graphicx}
\usepackage{xcolor}
\usepackage{placeins}
\usepackage{hyperref}
\usepackage{cleveref}
\usepackage{amsmath,mathtools}

\hypersetup{colorlinks=true, linkcolor=blue, citecolor=blue, urlcolor=blue}
\DeclareUrlCommand\path{\urlstyle{same}}

\newcommand{\rmnum}[1]{\romannumeral #1}
\DeclareMathOperator*{\argmin}{arg\,min}
\DeclareMathOperator*{\argmax}{arg\,max}
\DeclareMathOperator{\diag}{diag}
\DeclareMathOperator{\diver}{div}

\title{Random feature policy iteration for inverse mean field games with complex geometry and high dimensions}
\author{Xiangjian Zhang\thanks{Department of Mathematics,
University of Pennsylvania, Philadelphia, PA 19104; \texttt{chkirin@sas.upenn.edu}}
\and Shanyin Tong\thanks{Department of Mathematics,
University of Pennsylvania, Philadelphia, PA 19104; \texttt{tong3@sas.upenn.edu}}
\and Xiaofeng Xu\thanks{Department of Mathematics,
University of Pennsylvania, Philadelphia, PA 19104; \texttt{xkx5060@sas.upenn.edu}}
}
\date{\today}

\begin{document}
\maketitle

\begin{abstract}
We propose a random feature policy iteration method for reconstructing the obstacle function in a mean-field game (MFG) from observations of its value function, leveraging the method's mesh-free nature and flexibility to handle complex geometries and high-dimensional problems.
The method builds on a policy iteration framework for inverse MFGs, which decouples the coupled nonlinear forward–backward PDE system into a sequence of linear PDE and inverse problem solves. The fixed random feature representations preserve this linear structure, allowing the resulting subproblems to be solved by minimizing residuals at sampled collocation points through linear least-squares problems. The mesh-free nature of the method further extends the applicability of the policy iteration framework to inverse MFGs with complex geometry and in high dimensions. Another key advantage is its flexibility in using different spatial approximations for different variables: the obstacle function, value function, and population density can be represented in separate feature spaces, accommodating obstacles that are substantially less regular than the corresponding PDE solutions. Moreover, boundary conditions can be enforced through collocation on irregular domains without constructing boundary-fitted meshes. Numerical experiments involving smooth and discontinuous obstacles, a domain with holes, and higher-dimensional problems demonstrate the accuracy and flexibility of the proposed method.
\end{abstract}


\section{Introduction}
\label{sec:intro}
Mean-field games (MFGs) model strategic interactions in large populations of agents, where an individual's influence is negligible but the collective behavior of the population significantly affects each agent's decisions~\cite{LasryLions2007}. Each agent chooses a strategy to minimize a cost that typically depends on its state, control, and the population distribution. At equilibrium, the value function, representing an individual agent's optimal cost, satisfies a backward Hamilton--Jacobi--Bellman (HJB) equation coupled with a forward Fokker--Planck (FP) equation for the population density. Applications of MFGs include economics and finance~\cite{AchdouEtAl2022,CarmonaDelarueLacker2017} and traffic flow~\cite{HuangEtAl2021}; related mean-field control models have also been used for epidemic modeling~\cite{LeeEtAl2021}. The coupled PDE structure has motivated various numerical methods, including finite-difference schemes~\cite{AchdouCapuzzoDolcetta2010}, machine-learning-based approaches for higher-dimensional problems~\cite{RuthottoEtAl2020}.

In many applications, the cost function contains environmental information that is difficult to obtain directly. Inverse MFG problems seek to recover such unknown quantities from observations of equilibrium behavior, such as population densities, optimal strategies, or value functions. Recent studies have investigated uniqueness and stability~\cite{LiuMouZhang2023,RenEtAl2024}, as well as numerical reconstruction using PDE-constrained optimization~\cite{DingEtAl2022}, operator splitting~\cite{ChowEtAl2023}, bilevel optimization~\cite{YuEtAl2024}, and Gaussian processes~\cite{GuoEtAl2024}. In this paper, we consider the reconstruction of a spatially dependent obstacle function \(b(x)\), also called a potential, from linear measurements of the value function.
The obstacle enters the agents' running cost additively, making regions with larger values more costly to visit. Our numerical experiments use observations of the value function at the initial time.

The numerical reconstruction of an obstacle can be particularly challenging when its spatial structure is complicated, or the state space is high-dimensional or has complex geometry. Grid-based methods are effective on regular domains and in low dimensions, but constructing meshes for domains with irregular boundaries can be cumbersome, and the number of spatial degrees of freedom can grow rapidly with dimension.
Moreover, an obstacle with sharp variations or discontinuities may require a richer spatial representation than the value function or population density. The inverse policy iteration framework~\cite{RenSoedjakTong2026} provides linear PDE and inverse subproblems at each iteration, but its finite-difference implementation does not directly address these spatial discretization challenges.
These considerations motivate a mesh-free approximation that preserves the linear structure of the subproblems, allows different unknown functions to be represented using distinct feature functions, and accommodates irregular domains and high-dimensional settings.

Random features provide a natural way to construct such an approximation. Originally introduced for finite-dimensional approximations of kernel methods~\cite{RahimiRecht2007}, the basic idea is to represent an unknown function by a linear combination of randomly chosen functions whose parameters remain fixed after sampling. Only the linear coefficients of this combination need to be determined. For PDEs, random feature methods use collocation to fit the governing equations and boundary conditions at sampled points~\cite{ChenChiEYang2022,ChenELuo2023,ChiChenYang2024}. When the PDE operators are linear, the unknown coefficients can be determined by solving linear least-squares problems, without training the feature parameters. Since collocation does not require mesh connectivity, it offers flexibility in handling irregular geometries and higher-dimensional state spaces. For forward MFGs, Gaussian processes have been used to compute equilibrium solutions~\cite{MouYangZhou2022}, while parametrization method have been used to approximate nonlocal interaction kernels~\cite{AgrawalEtAl2022}. These approaches inspire us to use fixed features to discretize the state variables and reconstruct the unknown obstacle in an inverse MFG problem.

Motivated by these observations, we propose a random feature policy iteration method for inverse MFGs.
Building on~\cite{RenSoedjakTong2026}, we combine random feature approximation in space with implicit time discretization. The value function, population density, and obstacle are represented separately, and the discretized PDEs (linearized FP and HJB equations) are enforced by least-squares collocation. For a fixed policy, the resulting PDE problems and the inverse subproblem are solved using linear least-squares, without training the feature parameters or introducing nested nonlinear optimization. The feature values and their spatial derivatives are precomputed and reused throughout the iteration. For the periodic settings considered in our experiments, we construct the approximation spaces from randomly sampled Fourier modes, which directly satisfy the periodic boundary conditions. In this way, the method retains the linear structure of the subproblems while providing greater freedom in the choice of spatial representations.
Numerical experiments with smooth and discontinuous obstacles, a domain with holes, and problems in up to ten dimensions illustrate the method’s accuracy and flexibility.

The main \emph{contributions} of this paper are: (1) a mesh-free random feature discretization of inverse policy iteration that preserves linear least-squares PDE solves and the linear structure of the inverse step; (2) independently chosen feature spaces for the obstacle, value function, and population density, allowing the obstacle representation to be enriched without increasing the sizes of the other two spaces; and (3) a collocation framework for dealing with irregular boundaries and higher-dimensional state spaces without constructing a boundary-fitted spatial mesh. The central advantage is the flexibility to adapt the spatial representation to the unknown functions while retaining the computational simplicity of fixed-feature linear subproblems.

\section{Problem setup}
\label{sec:setup}
In this paper, we study an inverse mean field game (MFG) problem in which the goal is to recover a spatially dependent obstacle function \(b(x)\) from observations of the value function \(u\). Given a linear observation operator \(\mathcal G\) and spatial measurements \(g_i=g(x_i^{\rm obs})\) at the observation points \(\{x_i^{\rm obs}\}_{i=1}^{N_{\rm obs}}\),
we formulate the reconstruction as the following least-squares problem constrained by the MFG system:
\begin{subequations}\label{eq:opt}
\begin{alignat}{4}
&\underset{b(x)}{\operatorname{minimize}}
&&\quad \mathrlap{
    \mathcal J(b)
    :=\sum_{i=1}^{N_{\rm obs}}
    \left[\mathcal G u(x_i^{\rm obs})-g_i\right]^2}
&&
&&\label{eq:1a}\\
&\text{subject to}
&&\quad\smash{\raisebox{-1.5\baselineskip}{$
    \displaystyle\left\{
    \vphantom{\begin{aligned}
        \strut\\\strut\\\strut\\\strut
    \end{aligned}}
    \right.$}}
&&\;-\partial_tu-\varepsilon\Delta u
    +H(\nabla u)=b(x)+F(m)
&&\qquad\text{in }Q, \label{eq:1b}\\
&
&&
&&\;\partial_tm-\varepsilon\Delta m
    -\operatorname{div}\bigl(mH_p(\nabla u)\bigr)=0
&&\qquad\text{in }Q, \label{eq:1c}\\
&
&&
&&\;u(x,T)=u_T(x),\quad m(x,0)=m_0(x)
&&\qquad\text{in }\Omega, \label{eq:1d}\\
&
&&
&&\;\text{boundary conditions for }u\text{ and }m
&&\qquad\text{on }\partial\Omega\times(0,T). \label{eq:1e}
\end{alignat}
\end{subequations}
Here \(\Omega\subset\mathbb R^d\) is a bounded domain, \(Q:=\Omega\times(0,T)\), and \(T>0\) is the final time.
A common setting of state space in the MFG is the flat torus \(\mathbb T^d:=\mathbb R^d/\mathbb Z^d\), which can be identified with the unit cube \([0,1]^d\) under periodic boundary conditions. In most of the following discussions, we adopt this setting, but we will address specific boundary conditions separately in certain cases.

The MFG system~\eqref{eq:1b}--\eqref{eq:1e} consists of two coupled nonlinear PDEs. The backward Hamilton--Jacobi--Bellman
(HJB) equation~\eqref{eq:1b} determines the value function \(u(x,t)\), which represents the minimum expected cost for an agent starting from state \(x\) at time \(t\). The forward Fokker--Planck (FP) equation~\eqref{eq:1c} describes the evolution of the population density \(m(x,t)\ge 0\). The diffusion coefficient \(\varepsilon>0\) accounts for random fluctuations.
The Hamiltonian \(H:\mathbb R^d\to\mathbb R\) describes the cost associated with agents' control decisions, with \(H_p\) denoting its gradient, applying to \(p=\nabla u\), while \(F\) models interactions through the population density. The terminal cost \(u_T\) and initial density \(m_0\) are prescribed, with \(\int_\Omega m_0(x)\,dx=1\). Throughout, \(H\), \(F\), \(u_T\), and \(m_0\) are assumed known, while \(b\) is to be reconstructed from observations satisfying \(g\approx\mathcal G u\). In our numerical experiments, we consider the initial-time observation \(\mathcal G u=u(\cdot,0)\).

The obstacle \(b\) represents a state-dependent cost, making regions where its value is large more costly for agents to visit. Its spatial structure can be complicated: \(b\) may vary smoothly, contain multiple localized peaks or wells, or have discontinuities along irregular interfaces. Moreover, applications may involve domains with complex boundaries or high-dimensional state spaces, for which grid-based discretizations can become difficult or prohibitively expensive. These challenges motivate us to  develop a flexible, mesh-free method that can accommodate complex geometry, allow \(b\), \(u\), and \(m\) to be approximated independently, and avoid constructing a full tensor-product grid in high-dimensional problems.

\section{Preliminaries}
\subsection{Policy iteration for inverse MFGs}
\label{sec:inverse-pi}
Policy iteration is a classical method for optimal control problems; see~\cite{Bellman1957,Bertsekas1987,Bertsekas2011,Farahmand2008,Howard1960,LagoudakisParr2003,PutermanBrumelle1979,SantosRust2004,Scherrer2014,SuttonBarto2018} and the references therein.
Its basic idea is to alternate between evaluating the value function for a fixed policy and improving the policy using the resulting value function. In MFGs, the policy also governs the evolution of the population density. Fixing the policy decouples the forward Fokker--Planck (FP) and backward Hamilton--Jacobi--Bellman (HJB) equations, which can then be solved in succession~\cite{Cacace2021}. For inverse MFGs, the value-function step is extended to reconstruct the unknown obstacle function from observations~\cite{RenSoedjakTong2026}.

Let \(q:Q\to\mathbb R^d\) denote a policy representing a candidate control. The Lagrangian \(L:\mathbb R^d\to\mathbb R\) associated with the Hamiltonian \(H\) is defined by the Legendre transform \(L(q):=\sup_{p\in\mathbb R^d}(q\cdot p-H(p))\), and the optimum is obtained when $q=H_p(\nabla u)$. For a fixed policy $q$, the transport velocity \(H_p(\nabla u)\) in the FP equation~\eqref{eq:1c} is replaced by \(q\), while the nonlinear Hamiltonian term \(H(\nabla u)\) in the HJB equation~\eqref{eq:1b} is replaced by \(q\cdot\nabla u-L(q)\). This decouples the two PDEs in MFGs and leads to the following policy iteration for inverse MFGs~\cite{RenSoedjakTong2026}.

Starting from an initial policy \(q^{(0)}\in C^{1,0}(Q;\mathbb R^d)\), the inverse policy iteration consists of the following three steps for each \(k\ge0\):
\begin{enumerate}
\item[(\rmnum{1})] For the current policy \(q^{(k)}\), compute the density \(m^{(k)}\) by solving the linear FP equation:
\begin{equation}
\label{eq:pi-fp}
\begin{cases}
\partial_tm^{(k)}-\varepsilon\Delta m^{(k)}
-\diver(m^{(k)}q^{(k)})=0 &\text{in }Q,\\
m^{(k)}(x,0)=m_0(x) &\text{on }\mathbb T^d.
\end{cases}
\end{equation}
\item[(\rmnum{2})] Given \(q^{(k)}\) and \(m^{(k)}\), determine the obstacle function \(b^{(k)}\) and the corresponding value function \(u^{(k)}\) from observations \(\{(x_i^{\rm obs},g_i)\}_{i=1}^{N_{\rm obs}}\) by minimizing the data mismatch subject to a linear HJB equation:
\begin{subequations}\label{eq:linear-HJB}
\begin{alignat}{2}
&\underset{\mathclap{b^{(k)}(x)}}{\operatorname{minimize}}
&&\quad
\mathcal J(b^{(k)}):=
\sum_{i=1}^{N_{\rm obs}}\left[\mathcal Gu^{(k)}(x_i^{\rm obs})-g_i\right]^2
\label{eq:linear-HJB-obj}\\
&\operatorname{subject\ to}
&&\quad
\begin{cases}
-\partial_tu^{(k)}-\varepsilon\Delta u^{(k)}
+q^{(k)}\cdot\nabla u^{(k)}-L(q^{(k)})
=b^{(k)}+F(m^{(k)})
    &\text{in }Q,\\
u^{(k)}(x,T)=u_T(x)
    &\text{on }\mathbb T^d.
\end{cases}
\label{eq:linear-HJB-constraint}
\end{alignat}
\end{subequations}
Since \(q^{(k)}\) and \(m^{(k)}\) are fixed, the HJB constraint~\eqref{eq:linear-HJB-constraint} is linear in \(u^{(k)}\) and \(b^{(k)}\).
Hence \(u^{(k)}\), and consequently its linear observation \(\mathcal G u^{(k)}\), depends affinely on \(b^{(k)}\). Thus, this is a linear inverse problem.
The data-misfit objective is therefore convex and quadratic in the obstacle function, and step~(\rmnum{2}) is a linear least-squares problem.
\item[(\rmnum{3})] Update the policy using the value function \(u^{(k)}\) obtained in step~(\rmnum{2}):
\begin{equation}
q^{(k+1)}(x,t)=\argmax_{q\in\mathbb R^d}\left[q\cdot\nabla u^{(k)}-L(q)\right]=H_p(\nabla u^{(k)}).
\end{equation}
\end{enumerate}
Introducing the policy function decouples 
the strongly coupled forward-backward MFG systems, at the same time, also providing linear structures of the decoupled sub-PDEs and sub-inverse problems.
For a quadratic Hamiltonian \(H(p)=\tfrac12|p|^2\) and terminal-time derivative observations \(\mathcal G u=\partial_t u(\cdot,T)\), the theoretical analysis~\cite{RenSoedjakTong2026} guarantees its linear convergence under suitable regularity assumptions and for sufficiently small time horizons.

\subsection{Random feature approximation}
\label{sec:rf-approximation}
The random feature (RF)
approach~\cite{RahimiRecht2007,HuangZhuSiew2006,EMaWu2020,LamperskiLekang2024} parametrizes a function as a linear combination of randomly sampled feature functions. The parameters defining these functions are sampled and then fixed, leaving only their linear coefficients to be determined.

On a bounded domain \(\Omega\subset\mathbb R^d\),  random features can be written as 
\begin{equation}
\phi_j(x)=\sigma(w_j\cdot x+\beta_j),\qquad j=1,\ldots,M,
\end{equation}
where \(\sigma\) is a smooth nonlinear activation function, such as \(\tanh\), \(\sin\), or \(\cos\), and \((w_j,\beta_j)\in\mathbb R^d\times\mathbb R\) are sampled from a prescribed distribution. We approximate a function \(u\) by a linear combination of random features, 
\begin{equation}
u_M(x;\theta)=\sum_{j=1}^{M}\theta_j\phi_j(x)
=\Phi(x)^\top\theta,\qquad \theta\in\mathbb R^M,
\end{equation}
where \(\Phi(x):=(\phi_1(x),\ldots,\phi_M(x))^\top\). Given observations \(\{(x_i,g_i)\}_{i=1}^{N_{\rm obs}}\), where \(g_i\) is the observed value of \(u\) at \(x_i\in\Omega\), the coefficients \(\theta\) can be determined by solving a linear least-squares problem to minimize the total data misfit between data $g_i$ and its evaluation of random features at $x_i$,
\begin{equation}
\label{eq:rf-data-fit}
\theta^*=\argmin_{\theta\in\mathbb R^M}
\sum_{i=1}^{N_{\rm obs}}
\bigl[\Phi(x_i)^\top\theta-g_i\bigr]^2.
\end{equation}

The same approach can be used to solve PDEs by fitting the random-feature approximation to the governing equations and boundary conditions at collocation points~\cite{ChenChiEYang2022,DongLi2021,ChenELuo2023,ChiChenYang2024}. Consider the linear boundary-value problem
\begin{equation}
\begin{cases}
\mathcal Lu(x)=f(x) & \text{in }\Omega,\\
\mathcal Bu(x)=h(x) & \text{on }\partial\Omega,
\end{cases}
\end{equation}
where \(\mathcal L\) and \(\mathcal B\) are linear differential and boundary operators, respectively. Given interior collocation points \(\{x_i\}_{i=1}^{N_c}\subset\Omega\) and boundary points \(\{x_j^{\rm bd}\}_{j=1}^{N_{\rm bd}}\subset\partial\Omega\), the coefficients \(\theta\) are determined by minimizing the PDE and boundary residuals:
\begin{equation}
\label{eq:rf-collocation}
\begin{aligned}
\theta^*=\argmin_{\theta\in\mathbb R^M}\Biggl\{
&\sum_{i=1}^{N_{c}}
\bigl[(\mathcal L\Phi)(x_i)^\top\theta-f(x_i)\bigr]^2+\lambda\sum_{j=1}^{N_{\rm bd}}
\bigl[(\mathcal B\Phi)(x_j^{\rm bd})^\top\theta-h(x_j^{\rm bd})\bigr]^2
\Biggr\},
\end{aligned}
\end{equation}
where \(\lambda>0\) balances the PDE and boundary residuals. Since the operators are linear and the features are fixed, this is again a linear least-squares problem in \(\theta\).

For the inverse MFG problem on \(\mathbb T^d\), we use smooth features that are periodic in the spatial
variables, so the periodic boundary conditions are satisfied by construction.
For problems with other boundary conditions, the corresponding residuals can be enforced in the least-squares objective, as in the second term of~\eqref{eq:rf-collocation}.

\section{Random feature policy iteration for inverse MFG problem}
In this section, we combine random feature approximation in space with an implicit time discretization to implement the inverse policy iteration introduced in~\Cref{sec:inverse-pi}. For a fixed policy, the FP and HJB subproblems are linear and can be discretized as linear least-squares problems. We then convert the inverse HJB step to a least-squares problem involving only the obstacle coefficients.

\subsection{Discretization of the random feature policy iteration}
In this subsection, we discretize the inverse MFG system using random features in space and a uniform grid in time. For each \(v\in\{u,m,b\}\), let \(\Phi_v=(\phi_1^v,\ldots,\phi_{M_v}^v)^\top\) denote a fixed vector of nonzero, smooth periodic features. On the time grid \(t_n=n\Delta t\), \(n=0,\ldots,N_t\), with \(\Delta t=T/N_t\), we approximate the unknown functions by
\[
u(x,t_n)=:u^n(x)\approx \Phi_u(x)^\top\theta_u^n,
\qquad
m(x,t_n)=:m^n(x)\approx \Phi_m(x)^\top\theta_m^n,
\qquad
b(x)\approx \Phi_b(x)^\top\theta_b.
\]
Here, \(\theta_u^n\in\mathbb R^{M_u}\) and \(\theta_m^n\in\mathbb R^{M_m}\) depend on time, whereas \(\theta_b\in\mathbb R^{M_b}\) is time-independent. The periodic features enforce the spatial boundary conditions.

Let \(\{x_i\}_{i=1}^{N_c}\) and \(\{x_i^{\rm obs}\}_{i=1}^{N_{\rm obs}}\) denote
sampled collocation and observation points, respectively.
Both sets of points are fixed throughout all time steps and policy iterations.
We consider initial-time observations \(g_i=g(x_i^{\rm obs})\) of \(\mathcal Gu=u(\cdot,0)\).
Define the feature matrices by \((\mathbf\Phi_v)_{ij}=\phi_j^v(x_i)\) and \((\mathbf\Phi_{\rm obs})_{ij}=\phi_j^u(x_i^{\rm obs})\).
We write \((\nabla\mathbf\Phi_v)_{ij}=\nabla\phi_j^v(x_i)\) for the matrix whose columns are the feature gradients, and \((\Delta\mathbf\Phi_v)_{ij}=\Delta\phi_j^v(x_i)\) for the vector of feature Laplacians.

We first describe the discretization at a fixed policy iteration \(k\), suppressing the superscript \((k)\) on the policy and solution coefficients; in particular, \(q^n=q^{(k)}(\cdot,t_n)\).

\paragraph{Discretization of the linearized FP equation.} 
We first apply the RF approach to the initial condition of the FP equation, by approximating the initial density using the least-squares formulation~\eqref{eq:rf-collocation} with the identity operator:
\begin{equation}
\label{eq:initial-projection}
\theta_{m}^0=\argmin_{\theta\in\mathbb R^{M_m}}\sum_{i=1}^{N_{c}}\bigl[\Phi_m(x_i)^\top\theta-m_0(x_i)\bigr]^2.
\end{equation}
Applying the implicit Euler scheme to the linearized FP equation~\eqref{eq:pi-fp} gives, for \(n=0,\ldots,N_t-1\),
\begin{equation}
\label{eq:time-dis-fp}
\frac{m^{n+1}-m^n}{\Delta t}-\varepsilon\Delta m^{n+1}-\diver\bigl(m^{n+1}q^{n+1}\bigr)=0.
\end{equation}
Define \(\mathcal L_{\rm FP}^{n+1}:=I-\Delta t[\varepsilon\Delta+\diver(q^{n+1}\cdot)]\), 
so that~\eqref{eq:time-dis-fp} can be written as \(\mathcal L_{\rm FP}^{n+1}m^{n+1}=m^n\). 

At each subsequent time level, 
we apply the PDE-residual least-squares formulation~\eqref{eq:rf-collocation} to the operator \(\mathcal L_{\rm FP}^{n+1}\) and determine the density coefficients \(\theta_m^{n+1}\) by minimizing the residual of~\eqref{eq:time-dis-fp}:
\begin{equation}
\label{eq:lls-fp}
\theta_m^{n+1}
=\argmin_{\theta\in\mathbb R^{M_m}}
\sum_{i=1}^{N_{c}}\Bigl[
(\mathcal L_{\rm FP}^{n+1}\Phi_m)(x_i)^\top\theta
-\Phi_m(x_i)^\top\theta_m^{n}\Bigr]^2,
\end{equation}
for \(n=0,\ldots,N_t-1\). The corresponding normal equations are
\begin{equation}
\label{eq:s4-fp-normal}
(A_m^{n+1})^\top A_m^{n+1}\theta_m^{n+1}
=(A_m^{n+1})^\top\mathbf\Phi_m\theta_m^n,
\end{equation}
where the FP collocation matrix is obtained by applying the operator to each feature:
\begin{equation}
\label{eq:fp-operator}
\begin{aligned}
(A_m^{n+1})_{ij}
:&=(\mathcal L_{\rm FP}^{n+1}\phi_j^m)(x_i)\\
&=\phi_j^m(x_i)-\Delta t\Bigl[\varepsilon\Delta\phi_j^m(x_i)
+q^{n+1}(x_i)\cdot\nabla\phi_j^m(x_i)+(\diver q^{n+1})(x_i)\phi_j^m(x_i)\Bigr],
\end{aligned}
\end{equation}
using \(\diver(vq)=q\cdot\nabla v+v\diver q\).
Thus, the FP equation is solved through a sequence of linear least-squares problems forward in time.

\paragraph{Discretization of the linearized  HJB equation.} 
As with the initial density in FP, we impose the terminal condition of the HJB equation by the least-squares formulation~\eqref{eq:rf-collocation}:
\begin{equation}
\label{eq:terminal-projection}
\theta_u^{N_t}=\argmin_{\theta\in\mathbb R^{M_u}}\sum_{i=1}^{N_{c}}\bigl[\Phi_u(x_i)^\top\theta-u_T(x_i)\bigr]^2.
\end{equation}
With the policy \(\{q^n\}\) and density \(\{m^n\}\) fixed, we discretize the linearized HJB equation~\eqref{eq:linear-HJB} backward from \(u^{N_t}=u_T\). For \(n=N_t-1,\dots,0\),
\begin{equation}\label{eq:HJB-time-discrete}
-\frac{u^{n+1}-u^n}{\Delta t}-\varepsilon\Delta u^n+q^n\cdot\nabla u^n-L\bigl(q^n\bigr)=b+F\bigl(m^n\bigr).
\end{equation}
The corresponding linear operator is \(\mathcal L_{\rm HJB}^n:=\mathcal I-\Delta t[\varepsilon\Delta-q^n\cdot\nabla]\), and the RHS is \(f^n:=L(q^n) + F(m^n)\). Then~\eqref{eq:HJB-time-discrete} becomes \(\mathcal L_{\rm HJB}^nu^n=u^{n+1}+\Delta t\left[b+f^n\right]\).

For a candidate obstacle coefficient vector \(\theta_b\), let \(\theta_u^n(\theta_b)\) denote the coefficients obtained from the backward HJB solve. Set
\(f_i^n:=L(q^n(x_i))+F(\Phi_m(x_i)^\top\theta_m^n)\).
Starting from \(\theta_u^{N_t}(\theta_b)=\theta_u^{N_t}\), we apply the PDE-residual least-squares formulation~\eqref{eq:rf-collocation} to the operator \(\mathcal L_{\rm HJB}^n\) and determine the density coefficients \(\theta_u^n\) by minimizing the residual of~\eqref{eq:HJB-time-discrete}:
\begin{equation}
\label{eq:ls-hjb}
\begin{aligned}
\theta_u^n(\theta_b)=\argmin_{\theta\in\mathbb R^{M_u}}
\sum_{i=1}^{N_{c}}\Bigl[&
(\mathcal L_{\rm HJB}^n\Phi_u)(x_i)^\top\theta
-\Phi_u(x_i)^\top\theta_u^{n+1}(\theta_b)-\Delta t\,\Phi_b(x_i)^\top\theta_b-\Delta t\,f_i^n\Bigr]^2,
\end{aligned}
\end{equation}
for \(n=N_t-1,\ldots,0\). 
The corresponding normal equations are
\begin{equation}
\label{eq:hjb-normal}
(A_u^n)^\top A_u^n\theta_u^n(\theta_b)
=(A_u^n)^\top\bigl[\mathbf\Phi_u\theta_u^{n+1}(\theta_b)
+\Delta t\,\mathbf\Phi_b\theta_b+\Delta t\,\mathbf f^n\bigr],
\end{equation}
where
\begin{align}
(A_u^n)_{ij}
&:= (\mathcal L_{\rm HJB}^n\phi_j^u)(x_i)
= \phi_j^u(x_i)
-\Delta t\Bigl[
\varepsilon\Delta\phi_j^u(x_i)
-q^n(x_i)\cdot\nabla\phi_j^u(x_i)
\Bigr],
\label{eq:Au}
\\
(\mathbf f^n)_i
&:= f_i^n.
\label{eq:fn}
\end{align}

For \(L^2\)-orthogonal features, the corresponding collocation matrices are expected to have full column rank under sufficiently rich random sampling at collocation points.
If the least-squares systems are ill-conditioned or rank-deficient, we can add 
a Tikhonov-type smoothness regularization term \(\|\nabla v\|_{L^2}^2\), for \(v\in\{u,m\}\), to the corresponding objectives.
Under the random-feature approximation, this term takes the form \(\|\nabla\Phi_v\theta_v\|_{L^2}^2\) and can be approximated using the precomputed feature gradients \(\nabla\mathbf\Phi_v\).
The corresponding normal equation matrix then becomes \((A_v^n)^\top A_v^n+\gamma_v R_v\), where \(R_v:=N_c^{-1}\sum_{i=1}^{N_c}\nabla\Phi_v(x_i)^\top\nabla\Phi_v(x_i)\) and \(\gamma_v\geq0\) is the regularization parameter.

\paragraph{The discrete inverse problem and policy iteration.} 
After the above discussions of how the RF methods solve the linearized FP and HJB by linear systems, we now can combine them together with the policy iteration framework discussed in \Cref{sec:inverse-pi} for inverse MFGs.
We now summarize the random-feature policy iteration for the inverse MFG problem. After a one-time initialization, each policy iteration consists of a forward FP solve, an inverse HJB solve, and a pointwise policy update.

\emph{Precomputation and initialization.} Evaluate \(\mathbf\Phi_v\) for \(v\in\{u,m,b\}\), \(\mathbf\Phi_{\rm obs}\), and the first- and second-order spatial derivatives
\(\nabla\mathbf\Phi_v, \Delta\mathbf\Phi_v\) at the
sampled collocation points. Compute the initial density and terminal value coefficients \(\theta_m^0\) and \(\theta_u^{N_t}\) once by solving the linear system~\eqref{eq:initial-projection} and~\eqref{eq:terminal-projection}.
Initialize a smooth periodic initial policy and evaluated at the sampled collocation points \(\{q^{(0),n}(x_i): n=0,\ldots,N_t,\, i=1, \ldots, N_c\}\),
for \(k=0, 1, \ldots\), repeat the following three steps:
\begin{enumerate}
\item[(\rmnum{1})] \emph{Forward FP solve.} Given \(\{q^{(k),n}\}\), set the initial density coefficients \(\theta_m^{(k),0}=\theta_m^0\). For \(n=0,\ldots,N_t-1\), solve the linear least-squares,
\begin{equation}
\label{eq:pi-fp-pointwise}
\theta_m^{(k),n+1}
=\argmin_{\theta\in\mathbb R^{M_m}}
\sum_{i=1}^{N_{c}}\Bigl[
(\mathcal L_{\rm FP}^{(k),n+1}\Phi_m)(x_i)^\top\theta
-\Phi_m(x_i)^\top\theta_m^{(k),n}\Bigr]^2.
\end{equation}
\item[(\rmnum{2})] \emph{Inverse HJB step.} With the policy \(\{q^{(k),n}\}\)
and density coefficients \(\{\theta_m^{(k),n}\}\)
fixed, determine the obstacle coefficients \(\theta_b^{(k)}\) from the following linear inverse problem,
\begin{subequations}\label{eq:optimization}
\begin{alignat}{2}
&\operatorname*{minimize}_{\theta_b^{(k)}} \qquad
&&\mathcal J(\theta_b^{(k)})
:= \sum_{i=1}^{N_{\mathrm{obs}}}
\left[
\Phi_u(x_i^{\mathrm{obs}})^\top
\theta_u^{(k),0}(\theta_b)-g_i
\right]^2
\label{eq:objective}
\\[6pt]
&\text{subject to}\qquad
&&\begin{cases}
\theta_u^{k,N_t}
= \theta_u^{N_t},
\\[6pt]
\begin{aligned}[t]
\theta_u^{(k),n}(\theta_b^{(k)})
&= \operatorname*{arg\,min}_{\theta\in\mathbb R^{M_u}}
\sum_{i=1}^{N_c}
\Bigl[
(\mathcal L_{\mathrm{HJB}}^{(k),n}\Phi_u)(x_i)^\top
\theta-\Phi_u(x_i)^\top
\theta_u^{(k),n+1}(\theta_b^{(k)})
\\
&\qquad
-\Delta t\,\Phi_b(x_i)^\top\theta_b^{(k)}
-\Delta t\,f_i^{(k),n}
\Bigr]^2,
\end{aligned}
\\[6pt]
\text{for } n=N_t-1,\ldots,0.
\end{cases}
\label{eq:constraints}
\end{alignat}
\end{subequations}
\item[(\rmnum{3})] \emph{Policy update.} Update the policy at the collocation points: 
\begin{equation}
\label{eq:pi-policy}
q^{(k+1),n}(x_i)=H_p\bigl(\nabla\Phi_u(x_i)\theta_u^{(k),n}\bigr),
\qquad n=0,\ldots,N_t.
\end{equation}
For \(H(p)=\tfrac12|p|^2\), this reduces to 
\(q^{(k+1),n}(x_i)=\nabla\Phi_u(x_i)\theta_u^{(k),n}\) with precomputed feature gradients \(\nabla\mathbf\Phi_u\).
\end{enumerate}

\subsection{Reduction of the inverse step}
In this subsection, we explain the details on how to convert the constrained inverse problem~\eqref{eq:optimization} to an unconstrained linear least-squares problem in \(\theta_b^{(k)}\).
A common reduced-space approach to PDE-constrained inverse problems is to eliminate the state variables through the governing equations, leaving an optimization problem in the unknown parameters alone.
This method has been applied to inverse problems with linear parabolic PDE constraints~\cite{AdavaniBiros2008}; related variable-elimination techniques have also been developed for inverse PDE problems using random-weight neural networks~\cite{DongWang2023}.

Fix a policy iteration \(k\) and suppress the superscript \((k)\). Since the policy \(\{q^n\}\)
and density coefficients \(\{\theta_m^n\}\)
are fixed, 
the discretized HJB equation~\eqref{eq:ls-hjb} is linear in the value-function coefficients \(\{\theta_u^n\}\), and the obstacle coefficients \(\theta_b\) enters as an additive source.
The HJB collocation matrices \(\{A_u^n\}\)~\eqref{eq:Au}
and the known source terms \(\{\mathbf f^n\}\)~\eqref{eq:fn}
are independent of \(\theta_b\).
Each solution therefore depends linearly on its right-hand side. Starting from the terminal coefficients \(\theta_u^{N_t}\), which are independent of \(\theta_b\), the HJB normal equation~\eqref{eq:hjb-normal} at \(n=N_t-1\) gives
\begin{equation}
(A_u^{N_t-1})^\top A_u^{N_t-1}\theta_u^{N_t-1}(\theta_b)
=(A_u^{N_t-1})^\top\Bigl[\bigl(\mathbf\Phi_u\theta_u^{N_t}+\Delta t\,\mathbf\Phi_b\bigr)\theta_b+\Delta t\,\mathbf f^{N_t-1}\Bigr].
\end{equation}
The coefficient matrix is fixed, while the right-hand side is affine in \(\theta_b\). Hence, the solution has the form \(\theta_u^{N_t-1}(\theta_b)=S^{N_t-1}\theta_b+r^{N_t-1}\). By mathematical induction applying to~\eqref{eq:hjb-normal}, we therefore obtain the affine dependence of the value-function coefficients \(\{\theta_u^n\}\) on the obstacle coefficients \(\theta_b\) at every time level:
\begin{equation}
\label{eq:affine-representation}
\theta_u^n(\theta_b)=S^n\theta_b+r^n,
\qquad S^n\in\mathbb R^{M_u\times M_b},\quad r^n\in\mathbb R^{M_u}.
\end{equation}
The matrices \(\{S^n\}\) and vectors \(\{r^n\}\) are obtained from the backward recursions
\begin{equation}
\label{eq:multiple-rhs}
\begin{aligned}
(A_u^n)^\top A_u^n S^n
&=(A_u^n)^\top\bigl(\mathbf\Phi_u S^{n+1}+\Delta t\,\mathbf\Phi_b\bigr),\\
(A_u^n)^\top A_u^n r^n
&=(A_u^n)^\top\bigl(\mathbf\Phi_u r^{n+1}+\Delta t\,\mathbf f^n\bigr),
\end{aligned}
\end{equation}
for \(n=N_t-1,\ldots,0\), starting from \(S^{N_t}=0\) and \(r^{N_t}=\theta_u^{N_t}\). These are linear systems and can be solved sequentially to obtain \(S^0\) and \(r^0\).

Using the observation matrix \(\Phi_u(x_i^{\rm obs})\) and the data \(g_i\),
substituting~\eqref{eq:affine-representation} into~\eqref{eq:optimization} reduces the inverse step to an unconstrained linear least-squares problem,
\begin{equation}
\label{eq:reduced-normal}
\theta_b=\argmin_{\theta\in\mathbb R^{M_b}}\sum_{i=1}^{N_{\rm obs}}\left[\Phi_u(x_i^{\rm obs})^\top(S^0\theta+r^0)-g_i\right]^2.
\end{equation}
After determining \(\theta_b\), we evaluate the value-function coefficients at all time levels from~\eqref{eq:affine-representation}. Hence, for a fixed policy and density, the inverse HJB step consists of a backward sequence of linear solves and one unconstrained least-squares solve for the obstacle coefficients.

Combining this reduction with the forward FP solve~\eqref{eq:pi-fp-pointwise} and the pointwise policy update~\eqref{eq:pi-policy} gives the complete random-feature policy iteration summarized in Algorithm~\ref{alg:rf-inverse-mfg}.

\begin{algorithm}[H]
\caption{Random feature policy iteration for inverse MFGs}
\label{alg:rf-inverse-mfg}
\small
\begin{algorithmic}[1]
\Statex \textbf{Input:} MFG parameters \(\varepsilon,T,H,F,m_0,u_T\); observation points \(\{x_i^{\rm obs}\}\) with data \(\{g_i\}\);
sampling distributions for random features and collocation points; number of features \(M_u,M_m,M_b\); number of collocation points \(N_{c}\); number of time steps \(N_t\); tolerance \(\tau>0\); maximum number of iterations \(I_{\max}\ge1\).
\State \textbf{Precomputation and initialization:} 
\State Sample the features \(\Phi_u,\Phi_m,\Phi_b\) and collocation points \(\{x_i\}_{i=1}^{N_{c}}\); fix them throughout the iterations.
\State Set \(\Delta t=T/N_t\); precompute \(\mathbf\Phi_u,\mathbf\Phi_m,\mathbf\Phi_b,\mathbf\Phi_{\rm obs}\) and the spatial derivatives
\(\nabla\mathbf\Phi_u,\nabla\mathbf\Phi_m,\Delta\mathbf\Phi_u,\Delta\mathbf\Phi_m\).
\State Compute \(\theta_m^0\) and \(\theta_u^{N_t}\) from the linear systems~\eqref{eq:initial-projection} and~\eqref{eq:terminal-projection};
initialize \(q^{(0),n}=0\) for \(n=0,\ldots,N_t\).
\For{\(k=0,\ldots,I_{\max}-1\)}
  \State \textbf{Forward FP solve:} Set \(\theta_m^{(k),0}=\theta_m^0\) and solve the linear systems~\eqref{eq:pi-fp-pointwise} forward for \(n=0,\ldots,N_t-1\).
  \State \parbox[t]{\dimexpr\linewidth-\algorithmicindent\relax}{\textbf{Inverse HJB solve:} Set \(S^{(k),N_t}=0\), \(r^{(k),N_t}=\theta_u^{N_t}\); compute \(S^{(k),n}\) and \(r^{(k),n}\) backward from~\eqref{eq:multiple-rhs} with the current policy 
  \(\{q^{(k),n}\}\) and density coefficients \(\{\theta_m^{(k),n}\}\); solve the linear system~\eqref{eq:reduced-normal} for \(\theta_b^{(k)}\); evaluate \(\theta_u^{(k),n}=S^{(k),n}\theta_b^{(k)}+r^{(k),n}\) for \(n=0,\ldots,N_t\).
  }
  \State \textbf{Policy update:} Update \(q^{(k+1),n}\) using~\eqref{eq:pi-policy}.
  \If{\(\max_{0\le n\le N_t}
        \left(\sum_{i=1}^{N_c}
        \bigl|q^{(k+1),n}(x_i)-q^{(k),n}(x_i)\bigr|^2\right)^{1/2}
        \le\tau\)}
    \State \textbf{break}
  \EndIf
\EndFor
\Statex \textbf{Output:} Reconstructed obstacle function \(b(x)=\Phi_b(x)^\top\theta_b^{(k)}\).
\end{algorithmic}
\end{algorithm}

For noisy or incomplete observations, the reduced least-squares system~\eqref{eq:reduced-normal} may be ill-conditioned or rank-deficient. To stabilize the reconstruction, we can add
a Tikhonov-type smoothness regularization term \(\|\nabla b\|_{L^2}^2\) to the corresponding objective.
Under the random-feature approximation, this term takes the form \(\|\nabla\Phi_b\theta_b\|_{L^2}^2\) and can be approximated using the precomputed feature gradients
\(\nabla\mathbf\Phi_b\).
The corresponding normal equation matrix then becomes \((\mathbf\Phi_{\rm obs}S^0)^\top\mathbf\Phi_{\rm obs}S^0+\gamma_b R_b\), where \(R_b:=N_c^{-1}\sum_{i=1}^{N_c}\nabla\Phi_b(x_i)^\top\nabla\Phi_b(x_i)\) and \(\gamma_b\geq0\) is the regularization parameter.

\section{Numerical experiments}
\label{sec:experiments}
In this section, we illustrate the proposed random feature policy iteration method for inverse MFGs using two-dimensional examples, including discontinuous obstacles and domains with holes, as well as examples in five and ten dimensions. We also examine how the reconstruction error depends on the numbers of features and collocation points.

For all experiments, we set the final time \(T=1\), the diffusion coefficient \(\varepsilon=0.3\). Throughout, the coupling cost is \(F(m)=m^2\) and the Hamiltonian is \(H(\nabla u)=\frac12|\nabla u|^2\), so that \(L(q)=\frac12|q|^2\) and \(H_p(\nabla u)=\nabla u\). The terminal condition is \(u_T(x)=-m_0(x)\). We reconstruct the obstacle function \(b(x)\) from observations of the value function at the initial time, \(\mathcal Gu=u(x,0)\). To generate the synthetic data, we solve the MFG system~\eqref{eq:1b}--\eqref{eq:1e}
with the true obstacle function and evaluate the resulting initial-time value function at the observation points.
In the noisy-data experiment, we add independent Gaussian noise to each observation: $g_i^\sigma=g_i+0.01\|u(\cdot,0)\|_{L^2}\xi_i$, with 
$\xi_i\overset{\mathrm{i.i.d.}}{\sim}\mathcal N(0,1)$.
Thus, the standard deviation of the noise at each observation point is set to \(1\%\) of the \(L^2\)-norm of the noiseless initial-time value function.

For the spatial approximation of \(u\), \(m\), and \(b\), we use Fourier features with integer frequencies, so that they directly satisfies periodic boundary conditions and can be directly applied to the torus without enforcing extra boundary conditions. Given distinct nonzero frequency vectors \(\kappa_1,\ldots,\kappa_M\in\mathbb Z^d\), we define the feature vector
\begin{equation}\label{eq:features}
\Phi(x)=\bigl(1,\ \sin(2\pi\kappa_1\cdot x),\ldots,\sin(2\pi\kappa_M\cdot x),\ \cos(2\pi\kappa_1\cdot x),\ldots,\cos(2\pi\kappa_M\cdot x)\bigr)^\top,
\end{equation}
Each frequency contributes a sine--cosine pair, in addition to the constant feature. Hence, \(M\) sampled frequencies 
yield \(2M+1\) scalar features. Throughout this section, \(M\), \(M_u\), \(M_m\), and \(M_b\) denote numbers of sampled frequencies rather than numbers of scalar features.

We select the frequencies from a bounded subset of the integer lattice. Since the frequencies \(\kappa\) and \(-\kappa\) generate the same sine--cosine pair up to a sign, we retain only one representative of each pair. Uniform sampling from the resulting set would assign equal probabilities to all candidate frequencies. Instead, we favor lower frequencies: for sufficiently regular periodic functions, Fourier coefficients tend to decay as the frequency increases, so low-frequency modes often capture much of the large-scale variation. The higher-frequency candidates remain available to represent finer details, including sharp variations. We use polynomially decaying weights, as motivated by Fourier descriptions of periodic Sobolev spaces~\cite{KriegVyb2023}; related polynomial spectral densities have also been used in random Fourier feature approximations~\cite{PottsWeidensager2025}. Specifically, we define
\begin{equation}
\label{eq:sampling}
\begin{gathered}
\Lambda_K:=\bigl\{\kappa\in\mathbb Z^d\setminus\{0\}:\ \|\kappa\|_\infty\le K,\ \text{the first nonzero component of $\kappa$ is positive}\bigr\},\\
\rho_{\mathrm{iso}}(\kappa):=\bigl(1+|\kappa|^2\bigr)^{-\beta},\qquad
\rho_{\mathrm{prod}}(\kappa):=\prod_{\ell=1}^d\bigl(1+\kappa_\ell^2\bigr)^{-\beta},\qquad \kappa\in\Lambda_K,
\end{gathered}
\end{equation}
where \(K\in\mathbb N\) bounds the frequency components and \(\beta>0\) controls the preference for low frequencies. We sample frequencies without replacement from \(\Lambda_K\), selecting each remaining candidate with probability proportional to the chosen weight. In two dimensions, we use the isotropic weight \(\rho_{\mathrm{iso}}\) with \(\beta=2\) and \(K=25\). In five and ten dimensions, we use the product weight \(\rho_{\mathrm{prod}}\) with \(\beta=4\) and \(K=12\). The collocation points are sampled independently and uniformly from the spatial domain.

\subsection{Reconstruction of a smooth obstacle function}
\label{subsec:rec-2d-simple}
\begin{figure}[hp]
\centering
\includegraphics[trim=27.846 2 2 2]{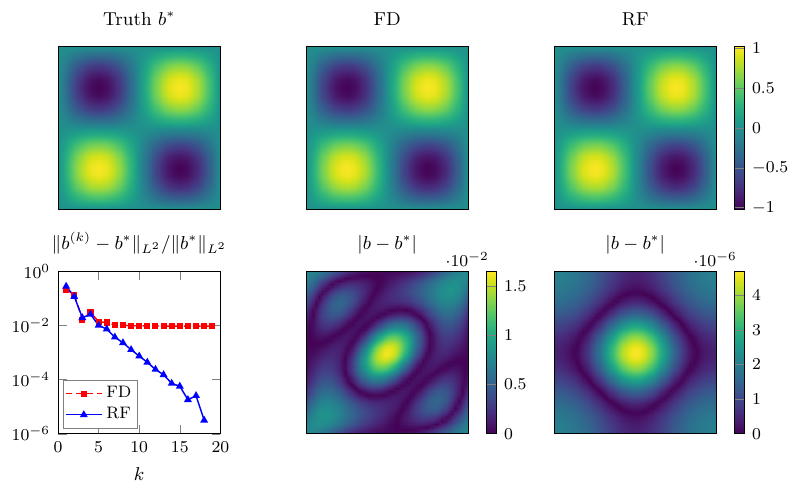}
\caption{Reconstruction results of the RF method and FD method for the two-dimensional inverse MFG problem~\eqref{eq:eq:2d-simple} with noiseless data. Top: reconstructed \(b\) for the two methods, compared with the true obstacle function \(b^*\). Bottom left: the decay of the relative error \(\|b^{(k)}-b^*\|_{L^2}/\|b^*\|_{L^2}\) with respect to the number of iterations \(k\) (displayed on a logarithmic scale on \(y\)-axis). Bottom middle and right: the pointwise absolute error \(\|b(x)-b^*(x)\|_{L^2}\) for the two methods.}
\label{fig:2d-simple-recover}

\vspace{1.2\baselineskip}
\includegraphics[trim=27.846 2 2 2]{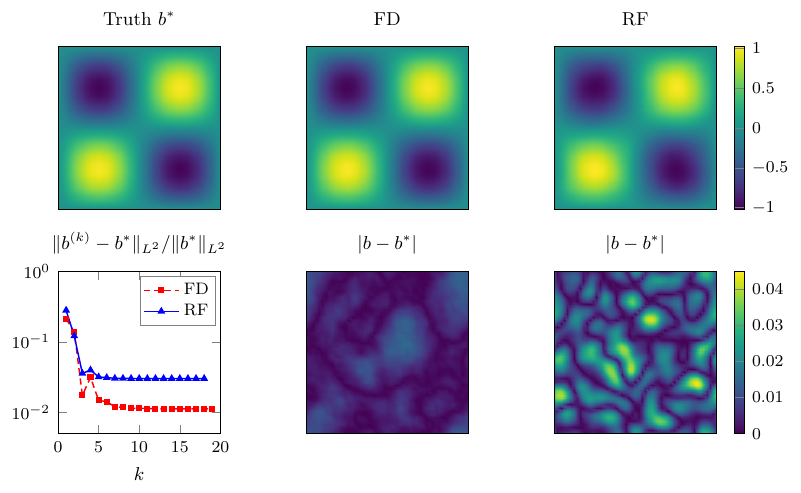}
\caption{Reconstruction results of the RF method and FD method for the two-dimensional inverse MFG problem~\eqref{eq:eq:2d-simple} with noisy data. Top: reconstructed \(b\) for the two methods, compared with the true obstacle function \(b^*\). Bottom left: the decay of the relative error \(\|b^{(k)}-b^*\|_{L^2}/\|b^*\|_{L^2}\) with respect to the number of iterations \(k\) (displayed on a logarithmic scale on \(y\)-axis). Bottom middle and right: the pointwise absolute error \(|b(x)-b^*(x)|\) for the two methods.
}
\label{fig:2d-simple-recover-noise}
\end{figure}
We first consider a smooth two-dimensional example as a basic validation of the proposed random feature (RF) policy iteration method. We compare it with the finite difference (FD) implementation of~\cite{RenSoedjakTong2026}. The true obstacle function is
\begin{equation}
\label{eq:eq:2d-simple}
b^*(x):=\sin(2\pi x_1)\sin(2\pi x_2),\qquad x\in\mathbb T^2,
\end{equation}
and the initial density is the following truncated periodization of a Gaussian,
\begin{equation}
\label{eq:2d-initial-dens}
m_0(x):=\frac{6.4}{\pi}\prod_{i=1}^{2}\sum_{j=-3}^{3}
\exp\!\bigl[-6.4\bigl(x_i-\tfrac12-j\bigr)^2\bigr],\qquad x\in\mathbb T^2,
\end{equation}
which satisfies \(\int_{\mathbb T^2}m_0(x)dx=1\).

We set the policy-iteration tolerance to \(\tau=10^{-8}\) and divide the time interval into \(N_t=100\) steps. The RF method uses \(N_c=2500\) collocation points and a common set of \(M=300\) sampled frequencies for \(u\), \(m\), and \(b\). Following~\cite{RenSoedjakTong2026}, the FD method uses the Engquist--Osher flux for the Hamiltonian and the transport term in the FP equation, with a uniform \(50\times50\) spatial grid. Both methods use observations of \(u(x,0)\) at the same \(2500\) grid nodes.

Figure~\ref{fig:2d-simple-recover} shows the reconstructions from noiseless data. 
Both methods capture the alternating peaks and troughs of the true obstacle, as shown in the top row of Figure~\ref{fig:2d-simple-recover}. However, the pointwise absolute error \(|b(x)-b^*(x)|\) maps in the bottom row show substantially smaller errors for RF, with the largest errors for both methods occurring near the center, where \(m_0\) and \(u_T\) are concentrated. The bottom left of Figure~\ref{fig:2d-simple-recover} shows the relative error \(\|b^{(k)}-b^*\|_{L^2}/\|b^*\|_{L^2}\) over the policy iterations in a logarithmic scale on the \(y\)-axis. The FD error levels off at around \(10^{-2}\), whereas the RF error continues to decrease to the order of \(10^{-6}\). 
This advantage is consistent with the low-frequency structure of the true obstacle: \(b^*\) is a linear combination of two Fourier modes with frequencies \((1,1)\) and \((1,-1)\), which are favored by our sampling scheme and can be represented accurately by Fourier features. In contrast, the FD scheme remains subject to spatial discretization and numerical-flux errors on the fixed grid.
\looseness=-1

We next examine the effect of observation noise. We perturb the observations of \(u(x,0)\) with independent Gaussian noise of standard deviation \(0.01\|u(\cdot,0)\|_{L^2}\) at each observation point. To stabilize the inverse problem in step~(\rmnum{2}), we add the Tikhonov regularization term \(\frac{\gamma}{2}\|\nabla b\|_{L^2}^2\) to the data misfit.
The regularization parameter is selected separately for each method using the L-curve criterion~\cite{Hansen1992,HansenOLeary1993}, which balances the fit to the observations against the spatial variation of the reconstructed obstacle. This gives \(\gamma=6.95\times10^{-8}\) for FD and \(\gamma=6.03\times10^{-8}\) for RF. All other settings remain unchanged.
As shown in Figure~\ref{fig:2d-simple-recover-noise}, both methods retain the main spatial structure of the obstacle.
However, the RF reconstruction exhibits more pronounced oscillatory errors, and its relative error levels off at around \(3\times10^{-2}\), compared with approximately \(10^{-2}\) for FD.
One possible explanation is that the Fourier-feature representation retains some of the oscillations introduced by the noisy data, whereas the
discretization of FD may help suppress them.
With the chosen regularization parameters, these noise-induced oscillations are not fully removed from the RF reconstruction.
Thus, although RF is more accurate for this smooth obstacle with noiseless observations, its advantage does not persist in this noisy-data experiment.
Traditional Tikhonov \(H_1\)-regularization may not be well suited to the RF setting, as the gradients of random features are not necessarily controlled by evaluating them at the collocation points.
This observation also motivates the development of more appropriate regularization strategies for RF methods, which we leave for future investigation.
\looseness=-1

\FloatBarrier
\subsection{Scaling law in two dimensions}
\label{subsec:2d-scaling}
We next investigate how the reconstruction error depends on the number of sampled frequencies \(M\) and the number of collocation points \(N_c\).
To examine the effect of enriching the feature space, we consider a smooth obstacle with a more complicated spatial structure:
\begin{equation}
\label{eq:2d-complex}
\begin{aligned}
b^*(x)&:=0.1\Bigl[\sin\bigl(2\pi x_1-0.6\sin(4\pi x_2)\bigr) +0.8\cos\bigl(2\pi x_2+0.5\sin(4\pi x_1)\bigr) \\
&\qquad+\exp\bigl(0.7\cos(2\pi x_1)+0.5\sin(2\pi x_2)+0.3\cos\bigl(2\pi(x_1-x_2)\bigr)\bigr) \\
&\qquad+0.35\sin\bigl(2\pi(x_1+x_2)\bigr)\cos\bigl(2\pi(2x_1-x_2)\bigr)\Bigr],\qquad x\in\mathbb T^{2},
\end{aligned}
\end{equation}
Unlike the obstacle in~\eqref{eq:eq:2d-simple}, which is a finite combination of a few Fourier modes, this function contains nonlinear compositions of trigonometric functions and an exponential term.
It therefore has a richer Fourier spectrum and cannot be represented exactly by a finite set of Fourier features.

\begin{figure}[htb]
\centering
\includegraphics[trim=2 2 2 2]{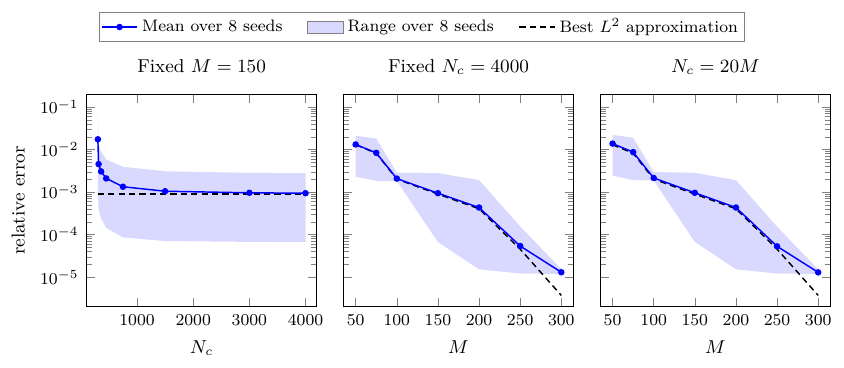}
\caption{Relative error \(\|b-b^*\|_{L^2}/\|b^*\|_{L^2}\) of the RF reconstruction with respect to different settings for the two-dimensional inverse MFG problem~\eqref{eq:2d-complex} (displayed on a logarithmic scale on \(y\)-axis). Solid lines: mean relative error over eight seeds. Shaded bands: range over the seeds. Dashed lines: relative error of the best \(L^2\) approximation of \(b^*\). Left: fix \(M=150\) and increase \(N_c\) from \(303\) to \(4000\). Middle: fix \(N_c=4000\) and increase \(M\) from \(50\) to \(3000\). Right: increase \(M\) from \(50\) to \(300\) together with \(N_c=20M\).}
\label{fig:2d-scaling-law}
\end{figure}

We consider three settings: (\rmnum{1}) fix \(M=150\) and increase \(N_c\) from \(303\) to \(4000\); (\rmnum{2}) fix \(N_c=4000\) and increase \(M\) from \(50\) to \(300\); and (\rmnum{3}) increase \(M\) from \(50\) to \(300\) with \(N_c=20M\), i.e., the fixed ratio between \(M\) and \(N_c\).
All other parameters are the same as in the noiseless experiment in Section~\ref{subsec:rec-2d-simple}.
We repeat each configuration with eight random seeds.
Figure~\ref{fig:2d-scaling-law} shows the mean relative reconstruction error and its range across seeds.
For comparison, the dashed curves show the relative error of the best \(L^2\) approximation of \(b^*\) in the sampled feature space for the obstacle.
These curves provide a lower bound on the reconstruction error attainable with the selected features.

When \(M=150\) is fixed, the reconstruction error decreases rapidly as \(N_c\) increases, but then levels off near \(10^{-3}\), close to the best-approximation error.
This suggests that once sufficiently many collocation points are used, the error is mainly limited by the fixed feature space, so adding more points brings little improvement.
In contrast, when \(N_c=4000\) is fixed, increasing \(M\) from \(50\) to \(300\) reduces the mean relative error by approximately three orders of magnitude, from \(10^{-2}\) to \(10^{-5}\).
The decrease is roughly exponential over the tested range, as indicated by the approximately linear trend on the logarithmic error scale, and closely follows the best-approximation error.
When \(N_c=20M\), the reconstruction error shows nearly the same decay rate and reaches comparable accuracy. Overall, the results illustrate that collocation points improve the discretization of a given feature space, while enriching the feature space is essential for further reducing the reconstruction error once the collocation error has become small.

\FloatBarrier
\subsection{Reconstruction of a discontinuous obstacle function}
\label{subsec:rec-2d-overlap}
We next consider a discontinuous obstacle to examine the recovery of sharp spatial transitions and, more importantly, the flexibility of choosing separate feature spaces for the value function, density, and obstacle.
The true obstacle \(b^*\) is piecewise polynomial on \(\mathbb T^2\), consisting of a constant background and three polynomial profiles supported on overlapping regions with irregular boundaries.
The contributions add where the regions overlap, producing jumps across their boundaries.
The complete definition is given in Appendix~\ref{app:discontinuous-obstacle}.

Since the obstacle contains discontinuities, it may require a richer Fourier representation than the value function and density.
We examine this possibility by keeping the number of sampled frequencies \(M_u=M_m=300\) fixed for $u$ and $m$ and comparing two obstacle feature spaces: (\rmnum{1}) \(M_b=300\) and (\rmnum{2}) a richer feature space for $b$ with \(M_b=1200\) sampled frequencies.
Both settings use \(N_c=10^4\) collocation points, and the remaining parameters follow the noiseless setup in Section~\ref{subsec:rec-2d-simple}.
To reduce oscillations while allowing sharp transitions, we add a total variation penalty \(\frac{\gamma}{2}\|\nabla b\|_{L^1}\) to the data-misfit objective in the inverse step, where \(\|\nabla b\|_{L^1}=\int_{\mathbb T^2}|\nabla b|\,dx\) is evaluated on a uniform \(100\times100\) spatial grid.
We use \(\gamma=1.78\times10^{-11}\) for \(M_b=300\) and \(\gamma=10^{-11}\) for \(M_b=1200\).
\begin{figure}[htb]
\centering
\includegraphics[trim=2 2 2 2]{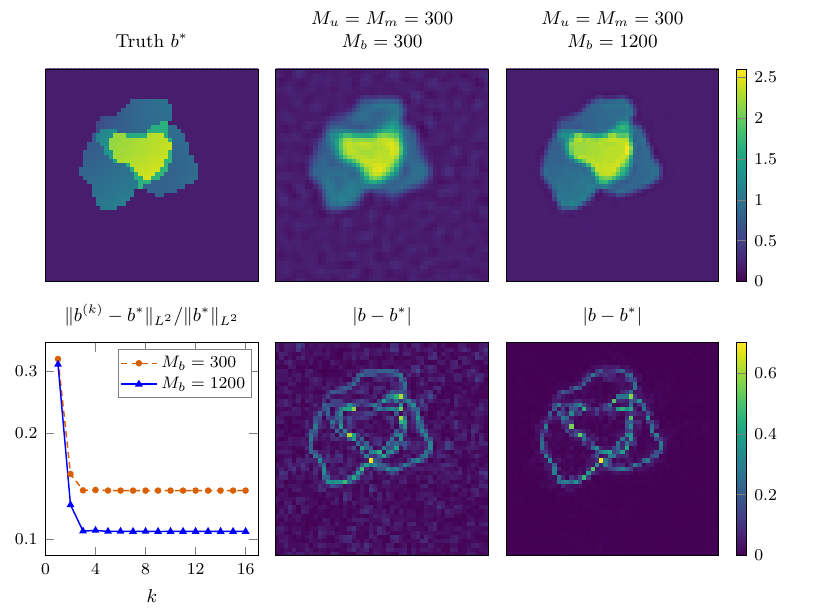}
\caption{Reconstruction results of the RF method in setting (\rmnum{1}), with \(M_u=M_m=M_b=300\), and setting (\rmnum{2}), with \(M_u=M_m=300\), \(M_b=1200\) for the two-dimensional inverse MFG problem~\eqref{eq:2d-overlap}. Top: reconstructed \(b\) for the two settings, compared with the true obstacle function \(b^*\). Bottom left: the decay of the relative error \(\|b^{(k)}-b^*\|_{L^2}/\|b^*\|_{L^2}\) with respect to the number of iterations \(k\) (displayed on a logarithmic scale on \(y\)-axis). Bottom middle and right: the pointwise absolute error \(|b(x)-b^*(x)|\) for the two settings.}
\label{fig:2d-overlap-recover}
\end{figure}

Figure~\ref{fig:2d-overlap-recover} shows that both settings recover the overall shape of the obstacle and the higher-valued region formed by the overlaps.
With \(M_b=1200\), the reconstructed boundaries are slightly sharper and the background is more uniform than with \(M_b=300\).
The relative reconstruction error decreases rapidly during the first few policy iterations and then levels off at approximately \(0.14\) for \(M_b=300\) and \(0.11\) for \(M_b=1200\).
The error maps show that the most prominent remaining errors occur near the discontinuities, which cannot be represented exactly by a finite Fourier expansion.
This comparison illustrates the advantage of assigning separate feature spaces to the MFG variables: the obstacle approximation can be enriched to better resolve its spatial structure without increasing the numbers of features used for \(u\) and \(m\).
If we have more knowledge about the function we want to reconstruct, we can even design or choose specific RF basis to better parametrize them.

\FloatBarrier
\subsection{Reconstruction on a torus with holes}
\label{subsec:rec-2d-holes}
\begin{figure}[htb]
\centering
\includegraphics[trim=2 2 2 2]{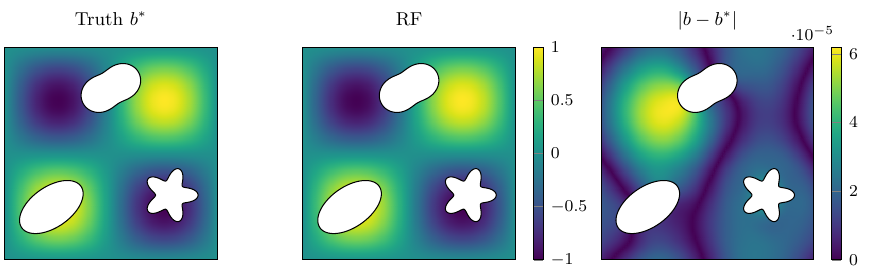}
\caption{Reconstruction result of the RF method for the two-dimensional inverse MFG problem~\eqref{eq:holes}. The holes are shown in white. Left: the true obstacle function \(b^*\). Middle: reconstructed \(b\) for the RF method. Right: the absolute error \(|b(x)-b^*(x)|\) for the RF method.}
\label{fig:2d-hole-recover}

\vspace{1.2\baselineskip}
\includegraphics[trim=2 2 2 2]{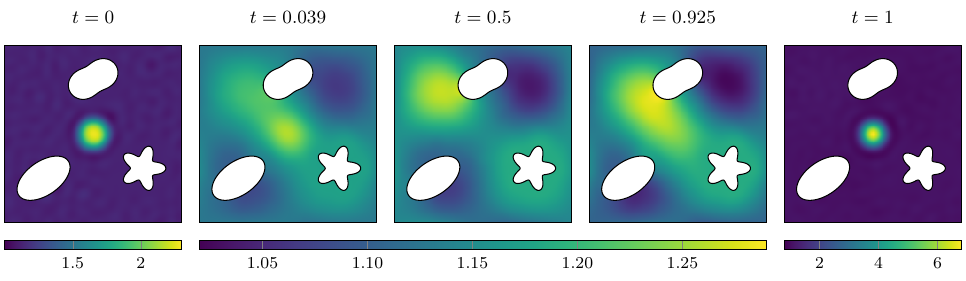}
\caption{Density evolution of the RF method for the two-dimensional inverse MFG problem~\eqref{eq:holes}, at \(t=0,0.039,0.5,0.925,1\) (left to right). The holes are shown in white. The three middle panels share one color scale; the panels at \(t=0\) and \(t=1\) have their own.}
\label{fig:2d-hole-density}
\end{figure}
We next test the RF method on a domain with three holes to examine its ability to handle irregular geometries.
We use the same smooth obstacle as in~\eqref{eq:eq:2d-simple}, now restricted to
\begin{equation}
\label{eq:holes}
\Omega:=\mathbb T^2\setminus(H_1\cup H_2\cup H_3),
\end{equation}
where \(H_1\) is a rotated ellipse and \(H_2\) and \(H_3\) have smooth, noncircular boundaries.
Their geometric definitions are given in Appendix~\ref{app:hole-geometry}.
The initial density consists of a constant background and a smooth bump centered at \(c=(\tfrac12,\tfrac12)\):
\begin{equation}
\label{eq:holes-m0}
m_0(x):=\frac{1}{Z}
\begin{cases}
1+\exp\!\left(1-\dfrac{1}{1-(|x-c|_{\mathbb T}/0.12)^2}\right),
& |x-c|_{\mathbb T}<0.12,\\[6pt]
1, & |x-c|_{\mathbb T}\ge 0.12,
\end{cases}
\qquad x\in\Omega.
\end{equation}
Here \(|x-c|_{\mathbb T}\) denotes the distance on the torus, and \(Z\) is the normalizing constant chosen so that \(\int_\Omega m_0(x)\,dx=1\).

On the boundaries of the holes, we impose a homogeneous Neumann condition for the value function and a zero-flux condition for the density:
\begin{equation}
\partial_n u=0,\qquad \varepsilon\partial_n m+m q\cdot n=0,
\qquad \text{on }\partial\Omega,
\end{equation}
where \(n\) is the outward unit normal to \(\Omega\).
We enforce these conditions at $N_{\rm bd}=1500$ boundary points, with $500$ points on each hole at equally spaced angular positions about its center.
As in Section~\ref{subsec:rec-2d-simple}, we use $N_c=2500$ interior collocation points, sampled uniformly from $\Omega$.
The observations of $u(x,0)$ and the reconstruction errors are evaluated at the $2185$ nodes of a $50\times50$ uniform spatial grid that lie in $\Omega$.

Figure~\ref{fig:2d-hole-recover} compares the reconstructed obstacle with the truth.
The two agree closely throughout $\Omega$, with pointwise errors on the order of $10^{-5}$ and no pronounced increase along the hole boundaries.
Figure~\ref{fig:2d-hole-density} shows the density computed by the RF method at $t=0$, $0.039$, $T/2$, $0.925$, and $T$.
The density is initially concentrated near the center, becomes more spread out at intermediate times, and develops a pronounced peak near the center by the final time.
Throughout the evolution, the density remains confined to the domain, consistent with the imposed zero-flux boundary condition.
Together, these results illustrate that RF collocation can handle the reconstruction and the density evolution on a domain with irregular holes by imposing boundary conditions at sampled points, without requiring a boundary-fitted mesh.

\FloatBarrier
\subsection{Scaling laws in five and ten dimensions}
\label{subsec:hd-scaling}
Finally, we examine the performance of the RF policy iteration method in five and ten spatial dimensions and study how the reconstruction error depends on the numbers of sampled frequencies \(M\) and collocation points \(N_c\).
In five dimensions, we consider the true obstacle function
\begin{equation}
\label{eq:5d-obstacle}
\begin{aligned}
b^*(x)&:=0.1\Bigl[\sin\bigl(2\pi x_1-0.6\sin(4\pi x_2)\bigr) +0.8\cos\bigl(2\pi x_3+0.5\sin(4\pi x_4)\bigr) \\
&\qquad+\exp\bigl(0.7\cos(2\pi x_1)+0.5\sin(2\pi x_5)+0.3\cos\bigl(2\pi(x_2-x_4)\bigr)\bigr) \\
&\qquad+0.35\sin\bigl(2\pi(x_1+x_3)\bigr)\cos\bigl(2\pi(2x_2-x_5)\bigr)\Bigr],\qquad x\in\mathbb T^{5}.
\end{aligned}
\end{equation}
The initial density is a product of truncated periodized Gaussians centered at the midpoint of the unit cube: 
\begin{equation}
\label{eq:5d-initial-dens}
m_0(x):=\left(\frac{4}{\pi}\right)^{5/2}\prod_{i=1}^{5}\sum_{j=-3}^3
\exp\!\bigl[-4\bigl(x_i-\tfrac12-j\bigr)^2\bigr],\qquad x\in\mathbb T^{5}.
\end{equation}
For the ten-dimensional experiment, we use the same obstacle as in~\eqref{eq:5d-obstacle}, extended to \(\mathbb T^{10}\) by making it independent of \(x_6,\ldots,x_{10}\).
The initial density, however, depends on all ten coordinates:
\begin{equation}
m_0(x):=\left(\frac{3}{\pi}\right)^{5}\prod_{i=1}^{10}\sum_{j=-4}^4
\exp\!\bigl[-3\bigl(x_i-\tfrac12-j\bigr)^2\bigr],\qquad x\in\mathbb T^{10}.
\end{equation}

\begin{figure}[htb]
\centering
\includegraphics[trim=2 2 2 2]{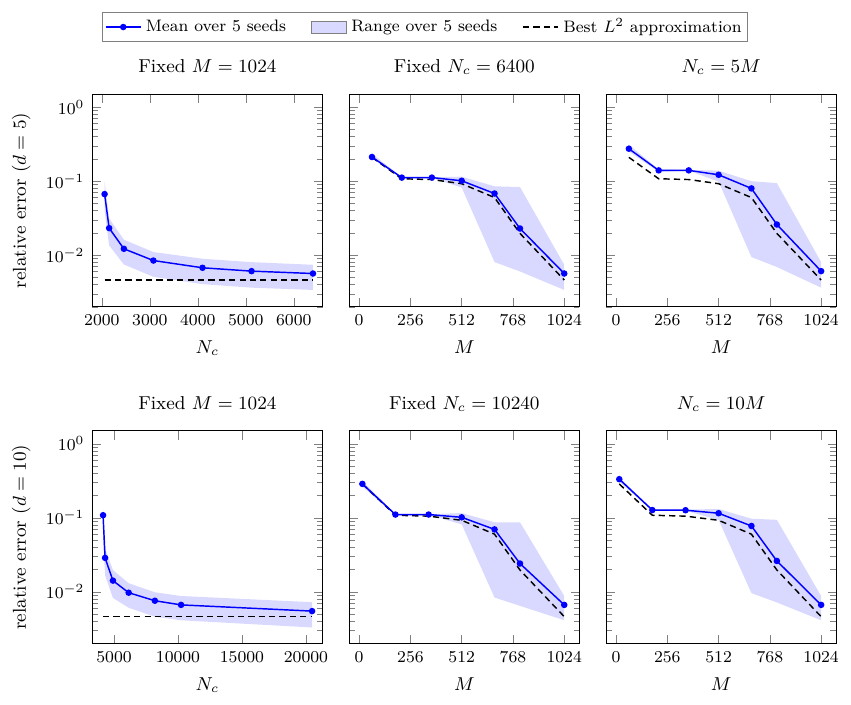}
\caption{Relative error \(\|b-b^*\|_{L^2}/\|b^*\|_{L^2}\) of the RF reconstruction with respect to difference settings for the five-dimensional (top) and ten-dimensional (bottom) inverse MFG problems~\eqref{eq:5d-obstacle} (displayed on a logarithmic scale on \(y\)-axis). Solid lines: mean relative error over five seeds. Shaded bands: range over the seeds. Dashed lines: relative error of the best \(L^2\) approximation of \(b^*\). Top left: fix \(M=1024\) and increase \(N_c\) from \(2057\) to \(6400\). Top middle: fix \(N_c=6400\) and increase \(M\) from \(64\) to \(1024\). Top right: increase \(M\) from \(64\) to \(1024\) together with \(N_c=5M\). Bottom left: fix \(M=1024\) and increase \(N_c\) from \(4145\) to \(20480\). Bottom middle: fix \(N_c=10240\) and increase \(M\) from \(16\) to \(1024\). Bottom right: increase \(M\) from \(16\) to \(1024\) together with \(N_c=10M\).}
\label{fig:hd-scaling-law}
\end{figure}

In both dimensions, observations are taken at 
\(N_{\rm obs}=16384\) uniformly sampled spatial points, and each configuration is repeated with five random seeds.
As in Section~\ref{subsec:2d-scaling}, we consider three settings: (\rmnum{1}) increase \(N_c\) while keeping \(M\) fixed; (\rmnum{2}) increase \(M\) while keeping \(N_c\) fixed; and (\rmnum{3}) increase \(M\) and \(N_c\) together at a fixed ratio.
For \(d=5\), we fix \(M=1024\) and vary \(N_c\) from \(2057\) to \(6400\) in setting~(\rmnum{1}); fix \(N_c=6400\) and vary \(M\) from \(64\) to \(1024\) in setting~(\rmnum{2}); and use \(N_c=5M\) and vary \(M\) from \(64\) to \(1024\) in setting~(\rmnum{3}).
For \(d=10\), the corresponding choices are \(M=1024\) with \(N_c\) varying from \(4145\) to \(20480\), \(N_c=10240\) with \(M\) varying from \(16\) to \(1024\), and \(N_c=10M\) with \(M\) varying from \(16\) to \(1024\).
Figure~\ref{fig:hd-scaling-law} shows the mean relative reconstruction error and its range across seeds, together with the best relative $L^2$ approximation error in the sampled obstacle feature space.

The top and bottom rows of Figure~\ref{fig:hd-scaling-law} show similar overall trends for the five- and ten-dimensional problems.
With \(M=1024\) fixed, increasing \(N_c\) initially reduces the reconstruction error substantially, while further improvements become small as the error approaches the best-approximation level.
When \(M\) increases, either at fixed \(N_c\) or with \(N_c\) proportional to \(M\), we observe that the error decreases relatively slowly for smaller \(M\) but begins to exhibit an exponential convergence trend at larger \(M\).
A similar phenomenon, where the asymptotic convergence rate becomes apparent only for sufficiently large numbers of features, has also been observed in high-dimensional function approximation using shallow ReLU$^k$ neural networks \cite{MaoXu2026}.
The two refinement settings give comparable accuracy overall and increasingly approach the best-approximation error as the feature space is enriched. 
The results demonstrate that RF policy iteration remains effective in these five- and ten-dimensional examples, while also illustrating the distinct roles of the two discretization parameters: additional collocation points improve the discretization within a fixed feature space, whereas additional frequencies enrich that space and lower its best-approximation error.

\section{Conclusion and discussion}
\label{sec:conclusion}
We develop a mesh-free random feature policy iteration method for reconstructing the obstacle function in inverse mean field games (MFGs) from observations of the value function.
Building on the inverse policy iteration framework, the method uses fixed random feature representations and collocation to preserve the linear structure of the decoupled PDEs and inverse subproblems.
For a fixed policy and density, the value-function coefficients depend affinely on the obstacle coefficients, allowing their elimination to reduce the inverse step to a linear least-squares problem.
Feature values and spatial derivatives are precomputed and reused throughout the iterations, and the algorithm combines linear solves with pointwise policy updates without training the random features.
A key advantage of the method is the use of separate feature spaces for the value function, population density, and obstacle, allowing the obstacle representation to be enriched independently of the other variables.
This flexibility improves the recovery of sharp interfaces in the discontinuous-obstacle experiment, although errors persist near the jumps.
Numerical experiments on a domain with holes and up to ten dimensions further illustrate the method’s ability to accommodate nontrivial geometries without boundary-fitted meshes and to extend to higher-dimensional settings while retaining the linear subproblem structure.

Several other directions are worth investigating. 
One is to select features adaptively from a larger candidate set rather than fix the active features in advance. Our numerical results indicate that the approximation capacity of the chosen feature spaces can limit reconstruction accuracy, motivating adaptive enrichment. The PDE residuals and observation mismatch could guide this process, identifying features that are most useful for reconstruction. Since finite Fourier representations cannot capture jump discontinuities exactly, incorporating features better suited to nonsmooth obstacles is also worth exploring.
Another direction is to develop regularization strategies tailored to random feature approximations.
Although we use Tikhonov regularization in the noisy-data experiment and total variation regularization for the discontinuous obstacle, oscillatory errors remain in the former, and errors persist near discontinuities in the latter. These results motivate the development of suitable penalties and parameter-selection strategies for random features in solving PDEs and inverse problems.

\section*{Acknowledgment}
ST is partially supported by National Science Foundation under Grant No.~DMS-2529292.

\appendix

\numberwithin{equation}{section}
\numberwithin{table}{section}
\section{Definitions for the numerical experiments}
\label{app:experiment-definitions}
In this section, we give the complete definitions of the true obstacle in Section~\ref{subsec:rec-2d-overlap} and of the domain in Section~\ref{subsec:rec-2d-holes}.

\subsection{Discontinuous obstacle function}
\label{app:discontinuous-obstacle}
For the discontinuous two-dimensional experiment, we construct an obstacle from three overlapping regions with irregular boundaries. The obstacle is defined by
\begin{equation}
\label{eq:2d-overlap}
b^*(x):=0.2+\sum_{i=1}^{3}\mathbf 1_{\Omega_i}(x)\,p_i\bigl(4(x-c_0)\bigr),
\qquad x\in\mathbb T^2,
\end{equation}
where $c_0=(0.45,0.60)$, $\mathbf{1}_{\Omega_i}$ is the indicator function of $\Omega_i$, and $p_i$ is a polynomial profile.
The constant $0.2$ gives the background value outside the three regions.
Within each region, the corresponding polynomial is added to the background, so the contributions accumulate where the regions overlap.

To define the regions $\Omega_i$, we start with rotated ellipses and perturb their boundaries in the radial direction:
\begin{equation}
\label{eq:overlap-regions}
\begin{split}
\Omega_i:=\bigl\{&c_i+R(\psi_i)\diag(\sigma_{i,1},\sigma_{i,2})\,
\varrho(\cos\vartheta,\sin\vartheta)^\top:\\
&\vartheta\in[0,2\pi),\quad 0\le\varrho<r_i(\vartheta)\bigr\},
\qquad i=1,2,3.
\end{split}
\end{equation}
Here $c_i$ is the center, $\sigma_{i,1}$ and $\sigma_{i,2}$ are the semi-axes of the underlying ellipse, and $\psi_i$ is its rotation angle.
The variables $\varrho$ and $\vartheta$ are radial and angular coordinates, respectively, and $r_i(\vartheta)$ describes the perturbation of the boundary.
The rotation matrix is given by
\begin{equation}
\label{eq:rotation}
R(\psi):=\begin{pmatrix}\cos\psi&-\sin\psi\\ \sin\psi&\cos\psi\end{pmatrix}.
\end{equation}
The geometric parameters are listed in Table~\ref{tab:overlap-regions}, and the radial functions are
\begin{equation}
\label{eq:overlap-radial}
\begin{aligned}
r_1(\vartheta)&:=1+0.16\cos(3\vartheta+0.25)+0.09\sin(5\vartheta-0.60),\\
r_2(\vartheta)&:=1+0.14\cos(4\vartheta-0.40)-0.10\sin(3\vartheta+0.80),\\
r_3(\vartheta)&:=1+0.15\sin(3\vartheta+0.70)+0.08\cos(6\vartheta-0.30).
\end{aligned}
\end{equation}
The resulting regions lie inside the unit square, overlap pairwise, and share a nonempty common intersection.
Since the regions are contained in the interior, the obstacle is compatible with periodicity on $\mathbb{T}^2$.
\begin{table}[h!]
\centering
\caption{Geometric parameters of the overlapping regions in~\eqref{eq:overlap-regions}.}
\label{tab:overlap-regions}
\begin{tabular}{@{}ccccc@{}}
\toprule
\(i\)&\(c_i\)&\(\sigma_{i,1}\)&\(\sigma_{i,2}\)&\(\psi_i\)\\
\midrule
\(1\)&\((0.36,0.54)\)&\(0.22\)&\(0.16\)&\(\pi/9\)\\
\(2\)&\((0.52,0.57)\)&\(0.20\)&\(0.15\)&\(-5\pi/36\)\\
\(3\)&\((0.43,0.68)\)&\(0.18\)&\(0.17\)&\(2\pi/9\)\\
\bottomrule
\end{tabular}
\end{table}
Using the common rescaled coordinates $\xi=(\xi_1,\xi_2)=4(x-c_0)$, the polynomial profiles are defined by
\begin{equation}
\label{eq:overlap-profiles}
\begin{aligned}
p_1(\xi)&:=0.75+0.18\,\xi_1-0.22\,\xi_2,\\
p_2(\xi)&:=0.65-0.12\,\xi_1+0.16\,\xi_2
+0.20\,\xi_1^2-0.10\,\xi_1\xi_2+0.14\,\xi_2^2,\\
p_3(\xi)&:=0.70+0.10\,\xi_1-0.08\,\xi_2
-0.16\,\xi_1^2+0.12\,\xi_1\xi_2+0.09\,\xi_2^2\\
&\quad+0.06\,\xi_1^3-0.04\,\xi_1^2\xi_2
+0.05\,\xi_1\xi_2^2-0.03\,\xi_2^3.
\end{aligned}
\end{equation}
Thus, the three regions carry linear, quadratic, and cubic polynomial profiles, respectively, and their boundaries introduce jump discontinuities in the obstacle.

\subsection{Domain with holes}
\label{app:hole-geometry}
The domain $\Omega$ in~\eqref{eq:holes} is obtained by removing three disjoint holes from $\mathbb T^2$. Their centers are
\[
c_1=(0.220,0.248),\qquad
c_2=(0.780,0.304),\qquad
c_3=(0.500,0.808).
\]
The first hole is a rotated ellipse with semi-axes $0.17$ and $0.095$:
\begin{equation}
\label{eq:hole-ellipse}
H_1:=\bigl\{c_1+R(7\pi/36)\diag(0.17,0.095)\,y:
 y\in\mathbb R^2,\ |y|\le1\bigr\},
\end{equation}
where $R$ is the rotation matrix defined in~\eqref{eq:rotation}.
The other two holes are defined by radial boundaries about their centers:
\begin{equation}
\label{eq:hole-radial-regions}
\begin{split}
H_i:=\bigl\{&c_i+\varrho(\cos\vartheta,\sin\vartheta)^\top: \vartheta\in[0,2\pi),\quad 0\le\varrho\le\rho_i(\vartheta)\bigr\},
\qquad i=2,3,
\end{split}
\end{equation}
where $\varrho$ is the radial coordinate and
\begin{equation}
\label{eq:hole-radial-functions}
\begin{aligned}
\rho_2(\vartheta)&:=0.100\bigl(1+0.3\cos(5\vartheta)\bigr),\\
\rho_3(\vartheta)&:=0.115\bigl(1+0.3\cos(2\vartheta-\pi/3)\bigr).
\end{aligned}
\end{equation}
The angular variations in $\rho_2$ and $\rho_3$ produce smooth, noncircular boundaries.
All three holes lie inside the unit square, so the outer edges retain the periodic identification of $\mathbb T^2$.

\bibliographystyle{elsarticle-num}
\bibliography{references}

\end{document}